\documentclass[12pt]{article}
\usepackage[english]{babel}
\usepackage{amsmath}
\usepackage[cp1251]{inputenc}
\usepackage{amscd}
\usepackage{amssymb}
\begin{document}

\noindent{\bf A Newman type bound for $L_p[-1,1]$-means
of the logarithmic derivative of polynomials having all
zeros on the unit circle}\\

{\bf Mikhail A. Komarov}\footnote{Vladimir State University, Gor'kogo street 87,
600000 Vladimir, Russia; \ \ e-mail: kami9@yandex.ru}\\

\noindent{\bf Abstract.} {\small
Let $g_n$, $n=1,2,\dots$, be the logarithmic derivative of a complex polynomial having all
zeros on the unit circle, i.e., 
a function of the form
$g_n(z)=(z-z_{1})^{-1}+\dots+(z-z_{n})^{-1}$, $|z_1|=\dots=|z_n|=1$.
For any $p>0$, we establish the bound
\[\int_{-1}^1 |g_n(x)|^p\, dx>C_p\, n^{p-1},\]
sharp in the order of the quantity $n$, where
$C_p>0$ is a constant, depending only on $p$. 
The particular case $p=1$ of this inequality
can be consi\-dered as a stronger variant of the well-known
estimate $\iint_{|z|<1} |g_n(z)|\,dxdy>c>0$
for the area integral of $g_n$, obtained by D.J. Newman (1972).
The result also shows that the set $\{g_n\}$
is not dense in the spaces $L_p[-1,1]$, $p\ge 1$.}\\

\noindent{\bf 1. Introduction. Main results}

\medskip

In 1971, C.K. Chui \cite{Chui71} asked whether there existed
an absolute constant $c>0$ such that
\begin{equation}\label{Newman}
    \iint_{|z|<1} |g_n(z)|\,dxdy> c
    \qquad (z=x+iy)
\end{equation}
for any rational function $g_n$, $n=1,2,\dots$, of the form
\begin{equation}\label{spf_N}
    g_n(z)=\sum_{k=1}^{n}\frac{1}{z-z_{k}}, \qquad |z_1|=\dots=|z_n|=1
\end{equation}
($g_n$ coincides with the logarithmic derivative of a polynomial
all of whose zeros lie on the unit circle).
The area integral in the left-hand side of (\ref{Newman})
shows the average strength in the unit disk
of the electrostatic field corresponding to
$n$ given point charges. 

In 1972, D.J. Newman \cite{Newman} has published
an elegant solution to Chui's problem
by showing that (\ref{Newman}) holds with
\[c:=\pi/18.\]
Later on, considering approximation properties of
logarithmic derivatives of complex po\-ly\-no\-mials in the Bers spaces
in bounded Jordan domains, Chui \cite{Chui} has proved,
in particular, that for every $q>2$ the set of sums (\ref{spf_N})
is dense in the space of functions $f$,
analytic in the open unit disk $D$, endowed with the norm
\[\|f\|_q:=\iint_{|z|<1} |f(z)|(1-|z|^2)^{q-2}\,dxdy.\]
So, Newman's bound (\ref{Newman}) shows 
that the condition $q>2$ cannot be improved here.

We also recall that it was conjectured \cite{Chui71}
that the integral in (\ref{Newman}) is minimized by the choice
$z_k=e^{2\pi k i/n}$ as the poles of $g_n$. This conjecture is still open,
however, its weighted $L_2$ analogs were recently solved in \cite{ABF}
as well as the problem of denseness of (\ref{spf_N})
in the corresponding function spaces in the unit disk.

Similar approximation problem for the unit interval
was recently posed in \cite[Sec. 4]{Bor-2016}:
whether sums (\ref{spf_N}) are dense in the complex space
$L_2[-1,1]$. The author \cite{K-Petrozav2019} has obtained the
negative answer to this question by showing the bound
\begin{equation}\label{Estimate-Petrozav}
    \int_{-1}^1 |g_n(x)|^2\,dx>\frac{1}{64}
\end{equation}
for any function $g_n$ of the form (\ref{spf_N}).

In this paper, using more subtle reasoning, we establish
a stronger result, showing, in particular, that $1/64$ in
(\ref{Estimate-Petrozav}) can be replaced by $n/800$.

\medskip

{\bf Theorem 1.} {\it For any function $g_n$ of the form $(\ref{spf_N})$
and any $p>0$ we have
\begin{equation}\label{th1}
    \int_{-1}^1 |g_n(x)|^p\, dx>\int_{-1}^1 |g_n(x)|^p\, |x|^p\, dx>C_p\, n^{p-1},
\end{equation}
where}
\[C_p:=\frac{3p^p(p+1)^{1-p}}{2^{p+5}(1+2p)^2}.\]

Choosing $p=1$ in (\ref{th1}),
we obtain the Newman type bound for 
the interval:
\begin{equation}\label{Newman type [-1,1]}
    \int_{-1}^1 |g_n(x)|\, dx>\int_{-1}^1 |xg_n(x)|\, dx>\frac{1}{192}
\end{equation}
and conclude that {\it rational functions $(\ref{spf_N})$
are not dense even in the space $L_1[-1,1]$}.


Note that the property (\ref{Newman type [-1,1]})
is stronger than (\ref{Newman}) in the sense that
Newman's result (\ref{Newman}) (although with a constant $c<\pi/18$)
immediately follows from (\ref{Newman type [-1,1]}):
\[\iint_{|z|<1} |g_n(z)|\,dxdy=\int_0^{2\pi}dt\int_0^1 |g_n(re^{it})|rdr
=\int_0^{\pi}dt\int_{-1}^1 |r g_n(re^{it})|dr>\frac{\pi}{192}.\]

The proof of Theorem 1 uses the following metric property
of functions (\ref{spf_N}) that is also of independent interest.
By $\mu(E)$ we will denote the linear measure of a set $E\subset \mathbb{R}$.

\bigskip

{\bf Theorem 2.} {\it For any $g_n$ of the form $(\ref{spf_N})$,
the measure $\mu(E)$ of the set
\[E=E_\delta(g_n)=\big\{x\in [-1,1]:\ \big|{\rm Re}(xg_n(x))\big|\ge \delta n\big\},
\qquad 0<\delta<1/2,\]
admits the estimate
\begin{equation}\label{th2}
    \mu(E)\ge K n^{-1}, \qquad
    K=K(\delta):=\frac{3}{32}\frac{1-2\delta}{(1+2\delta)^2}>0.
\end{equation}
Moreover, the set $E$ is concentrated near
the endpoints of the segment $[-1,1]$ in the sense that
$\mu(E)$ in $(\ref{th2})$ can be replaced by $\mu(E\cap \Delta)$, where}
\[\Delta=\Delta_{n,\delta}:=\Big\{x\in [-1,1]:\ |x|>1-\frac{3}{(2+4\delta)n}\Big\}.\]

In particular, by the definition of the set $E$ with $\delta:=p/(2(p+1))$, we have
\[\int_{-1}^1 |xg_n(x)|^p\, dx>\int_{E} \big|{\rm Re}(xg_n(x))\big|^p\,dx
\ge (\delta n)^{p}\mu(E)\ge \delta^pK(\delta)n^{p-1}=C_p\,n^{p-1},\]
and Theorem 1 follows.

\smallskip

Both estimates (\ref{th1}), (\ref{th2}) are sharp
in the order of the quantity $n$. Indeed, firstly, for the function
$g_n(x)=\tilde{g}_n(x):=nx^{n-1}/(x^n+i)$ we have
\[\int_{-1}^1 |\tilde{g}_n(x)|^p\, dx=
2n^{p-1}\int_0^1 \frac{(t^{1-\frac{1}{n}})^{p-1}dt}{(t^{2}+1)^{p/2}}
\le \tilde{C}_p\, n^{p-1}, \qquad \tilde{C}_p=\int_0^1 \frac{2t^{\kappa}\,dt}{(t^{2}+1)^{p/2}}\]
($p>0$, $\kappa=\min\{p-1;0\}$).
Secondly, if $0<\delta<1/2$ then from
\[{\rm Re}(x\tilde{g}_n(x))=\frac{nx^{2n}}{x^{2n}+1},
\qquad E_\delta(\tilde{g}_n)=
\big\{x\in [-1,1]:\ |x|\ge (\delta^{-1}-1)^{-1/(2n)}\big\}\]
and $a^{-t}>1-t\ln a$ ($a>1$) we deduce
\[\mu(E_\delta(\tilde{g}_n))<\tilde{K}n^{-1}, \qquad
\tilde{K}=\ln (\delta^{-1}-1)>0.\]

{\bf Remark 1.} For all $\delta\ge 1/2$, we have
\[\mu(E_\delta(\tilde{g}_n))=
\mu\big(\big\{x\in [-1,1]:\ x^{2n}\ge
\delta/(1-\delta)\big\}\big)=0,\]
therefore the property $\lim_{\delta\to 1/2} K(\delta)=0$
(see (\ref{th2})) is natural.

\smallskip

{\bf Remark 2.} 
Estimates of logarithmic derivatives of
polynomials with restricted zeros are well-known.
Thus, it was proved in \cite{Govorov-Lapenko}
that for polynomials $p_n$, whose zeros $z_1,\dots,z_n$
lie in the half-disk $U=\{z:|z|\le 1,\,{\rm Im}\,z\ge 0\}$,
the inequality $|p_n'(x)/p_n(x)|>\delta n$ ($\delta>0$)
holds for all $x\in [-1,1]$ outside an exceptional set
$e_\delta\subset [-1,1]$ of a measure
\begin{equation}\label{GovorovLapenko}
    \mu(e_\delta)<70e\delta.
\end{equation}
If, moreover, all $z_1,\dots,z_n\in [-1,1]$ then this 
estimate of $p_n'/p_n$ remains true with
\begin{equation}\label{GovorovLapenko[-1,1]}
    \mu(e_\delta)\le (\sqrt{1+4\delta^2}-1)/\delta
\end{equation}
(see \cite{Govorov-Lapenko}), 
and also the following weighted estimate holds (see \cite{K-Trudy}):
\begin{equation}\label{Trudy}
    \mu\big(\big\{x\in[-1,1]: \ (1-x^2)|p_n'(x)|
    \le \delta n |p_n(x)|\big\}\big)\le 2\delta,
\end{equation}
whose proof uses P. Borwein's result \cite{Borwein}.

Estimates (\ref{GovorovLapenko})--(\ref{Trudy})
are useful in constructing P. Tur\'{a}n type 
inequalities on a segment,
see \cite{K-AnalysisMath}, \cite{Erdelyi2020}, \cite{K-AnalysisMath21},
\cite{K-Trudy}. For example, it was shown by the author \cite{K-AnalysisMath}
that (\ref{GovorovLapenko}) implies
\begin{equation}\label{K-AnMath19}
    \|p_n'\|>C\sqrt{n}\,\|p_n\|, \qquad C:=2(3\sqrt{210e})^{-1},
    \qquad \|f\|:=\|f\|_{C[-1,1]},
\end{equation}
for polynomials $p_n$ having all $n$ zeros in the upper unit half-disk $U$.
T. Erd\'{e}lyi has considered $\mathcal{F}_{n,k}^c$
($k=0,1,\dots,n$), the classes of all complex polynomials
of degree at most $n$ having at least $n-k$ zeros in $U$,
and extended the result (\ref{K-AnMath19}) as follows \cite{Erdelyi2020}:

if $p_n\in \mathcal{F}_{n,k}^c$, $1\le k\le 10^{-5} n$, then we have
\begin{equation}\label{Erd n/100000}
    \|p_n'\|\ge \frac{1}{144e}\Big(\frac{n-k}{2k}\Big)^{1/2}\|p_n\|;
\end{equation}

for $p_n\in \mathcal{F}_{n,k}^c$, $1\le k\le n$, with at least one zero
in $[-1,1]$, we have
\[\|p_n'\|\ge \max\Big\{\frac{1}{2},\,\frac{1}{448}
\Big(\frac{n-k}{2k}\Big)^{1/2}\Big\}\|p_n\|.\]

{\bf Remark 3.} In connection with \cite{K-Petrozav2019},
V.I. Danchenko asked if (\ref{Estimate-Petrozav})
holds for sums (\ref{spf_N}) with $|z_k|\le 1$
rather than $|z_k|=1$ ($k=1,\dots,n$). It is natural to complete
this problem with the question on the existence
of an estimate of the form
\[\mu(G_\delta)\ge K_1(\delta)n^{-1}, \qquad
G_\delta:=\big\{x\in [-1,1]:\, |p_n'(x)|\ge \delta n |p_n(x)|\big\},
\quad \delta>0,\]
for the class $\mathcal{P}_n^c(\overline{D})$ of all
$n$th degree complex polynomials $p_n$,
all of whose zeros lie in the closed unit disk $|z|\le 1$.
In this issue, we only observe the fact that, simultaneously,
\begin{equation}\label{G delta [0,1]}
  \mu(G_\delta\cap [-1,0])>0 \quad \text{and} \quad \mu(G_\delta\cap [0,1])>0
\end{equation}
for any $\delta\in (0,1/2)$. Indeed, if $x=+1$ or $-1$ then
for $p_n\in \mathcal{P}_n^c(\overline{D})$
\[\left|\frac{p_n'(x)}{p_n(x)}\right|=\left|\sum_{k=1}^n \frac{1}{1-u_k}\right|
\ge \sum_{k=1}^n {\rm Re}\,\frac{1}{1-u_k}\ge \frac{n}{2}
\qquad (|u_k|\le 1),\]
where $n/2>\delta n$ because $\delta<1/2$.

\medskip

{\bf Corollary 1.} {\it If $p_n\in \mathcal{P}_n^c(\overline{D})$, then}
\[\|p_n'\|\ge \frac{1}{4}\|p_n\|.\]

{\it Proof.} Let $a\in [-1,1]$ be a point for which
$|p_n(a)|=\|p_n\|$. We will assume $a\in [0,1]$.

Put $G:=G_\delta$ with $\delta:=(2n)^{-1}$, i.e.,
$G=\{x\in [-1,1]: |p_n'(x)|\ge |p_n(x)|/2\}$.

If $a\in G$ then $\|p_n'\|\ge |p_n'(a)|\ge |p_n(a)|/2=\|p_n\|/2$.

If $a\in [0,1]\setminus G$ then, by (\ref{G delta [0,1]}),
there is a point $b\in G\cap [0,1]$ such that $|p_n'(b)|=|p_n(b)|/2$
and, for $x$ in the open interval between $a$ and $b$,
\[|p_n'(x)|<|p_n(x)|/2\le \|p_n\|/2; \qquad |b-a|\le 1.\]
From this, we have 
\[\|p_n'\|\ge |p_n'(b)|=\frac{|p_n(b)|}{2}\ge
\frac{1}{2}\Big(\|p_n\|-\left|\int_a^b p_n'(x)\,dx\right|\Big)\ge \frac{1}{4}\|p_n\|,\]
and Corollary 1 follows.

\smallskip


Now we can extend (\ref{K-AnMath19}) to the class of
polynomials with zeros in the unit disk:

\medskip

{\bf Corollary 2.} {\it If $p_n\in \mathcal{P}_n^c(\overline{D})$, then
\[\|p_n'\|\ge \max\Big\{\frac{1}{4},\,
\frac{1}{900}\sqrt{\frac{\max\{n^+,\,n^-\}+n_0}{2\min\{n^+,n^-\}+1}}\,\Big\}\|p_n\|,\]
where $n^+$ $(n^-, n_0)$ denotes a number of roots $z_k$
of $p_n$ whose imaginary part, ${\rm Im}\,z_k$, is
positive $($correspondingly, is negative, or zero$)$.}

\medskip

{\it Proof.} Because of symmetry and the result (\ref{K-AnMath19}), we may assume $n^+\ge n^-\ge 1$.
At the same time, $n_0+n^+=n-n^-$. Therefore,
the desired estimate takes the form
\[\|p_n'\|\ge \max\Big\{\frac{1}{4},\,M\Big\}\|p_n\|,
\qquad M:=\frac{1}{900}\Big(\frac{n-n^-}{1+2n^-}\Big)^{1/2}.\]
In the case $n^->10^{-5} n$, this estimate follows from Corollary 1
since in this case $M<1/4$. Now let $n^-\le 10^{-5} n$.
Obviously, $p_n\in \mathcal{F}_{n,k}^c$ with $k=n^-\ge 1$.
Thus, by (\ref{Erd n/100000}) we have
\[\|p_n'\|\ge \frac{1}{144e}
\Big(\frac{n-n^-}{2n^-}\Big)^{1/2}\|p_n\|>M\|p_n\|,\]
and the result follows from this and Corollary 1.

\bigskip

\noindent{\bf 2. Lemmas}

\medskip

For any $z_1,\dots,z_n$ on the unit circle $|z|=1$,
by using the identity
\[xg_n(x)=\sum_{k=1}^n \frac{x}{x-z_k}=-\frac{1}{2}\bigg(\sum_{k=1}^n
\frac{z_k+x}{z_k-x}-n\bigg),\]
we get
\begin{equation}\label{sum P - n}
    \big|{\rm Re}(xg_n(x))\big|
    =\frac{1}{2}\bigg|\sum_{k=1}^n P(z_k;x)-n\bigg|, \qquad -1\le x\le 1,
\end{equation}
where $P$ is the Poisson kernel
\[P(v;x):={\rm Re}\,\frac{v+x}{v-x}=
\frac{1-x^2}{1-2x\,{\rm Re}\,v+x^2}\ge 0, \qquad
-1\le x\le 1, \quad |v|=1.\]

It can be checked that the set of all points
$z\in \mathbb{C}$, satisfying the inequality ${\rm Re}\,((v+z)/(v-z))\ge h$,
coincides with the disk $|z-vh/(h+1)|\le 1/(h+1)$ for
any given $v$ ($|v|=1$) and $h>0$. We will be interested
in the intersection of such disks with the unit interval,
therefore, we investigate inequalities $P(v;x)\ge h$.

Fix $\rho\in (0,1/4]$ and put
\[T(h)=\sqrt{1+\rho^2-\frac{2\rho}{h}}, \qquad 1\le h\le (2\rho)^{-1}.\]
In the interval $1\le h\le (2\rho)^{-1}$,
the function $T(h)$ is well defined ($1+\rho^2-2\rho/h\ge (1-\rho)^2$)
and strictly increasing. In particular,
\[T(h)\le T((2\rho)^{-1})=\sqrt{1-3\rho^2}<1.\]

We formulate and prove two lemmas below for arbitrary
$\rho\in (0,1/4]$, however, in Sec. 3 we need only the case
$\rho=(4+8\delta)^{-1}n^{-1}$.

\medskip

{\bf Lemma 1.} {\it Let $|v|=1$, $\rho\in (0,1/4]$ and
\[{\rm Re}\,v\ge T=T(h) \quad \text{for some} \quad h\in [1,(2\rho)^{-1}].\]
Then the inequality
\[P(v;x)\ge h, \qquad x\in S,\]
holds, where $S=S(h)$ is the segment
\[S(h)=[x_-,x_+], \qquad
x_\pm:=\frac{\sqrt{h^2-2\rho h+\rho^2h^2}\pm (1-\rho h)}{h+1}.\]
At the same time, for any $h\in [1,(2\rho)^{-1}]$ the inclusion holds:}
\begin{equation}\label{S*}
    S(h)\supset S^*=\Big[\frac{\sqrt{1-3\rho^2}-\rho}{1+2\rho},\ 1-\rho\Big],
    \qquad \mu(S^*)>\frac{5\rho}{4}.
\end{equation}

{\bf Proof.} Fix $h\in [1,(2\rho)^{-1}]$.
The inequality $P(v;x)\ge h$ is equivalent to
\[(h+1)x^2-2hx\, {\rm Re}\,v+h-1\le 0,\]
or
\begin{equation}\label{(x-...)^2<}
    \Big(x-\frac{hd}{h+1}\Big)^2\le \frac{h^2(d^2-1)+1}{(h+1)^2}, \qquad d:={\rm Re}\,v.
\end{equation}
Here by $d\ge T=T(h)$ and $\rho h\le 1/2$ we have
\[h^2(d^2-1)+1\ge h^2(T^2-1)+1=(1-\rho h)^2>0.\]
Hence, for every individual value of $d$, $T\le d\le 1$,
the segment (\ref{(x-...)^2<}) is non-degenerate
and has the endpoints
\[A(d)=\frac{hd-\sqrt{h^2(d^2-1)+1}}{h+1}, \qquad
B(d)=\frac{hd+\sqrt{h^2(d^2-1)+1}}{h+1}.\]

Clear that at points
\[x\in \bigcap_{T\le d\le 1}[A(d),B(d)],\]
the inequality $P(v;x)\ge h$ holds for any $d={\rm Re}\,v\in [T,1]$.
But
\[\min\{B(d):T\le d\le 1\}=B(T)\]
and, because of the equality $A(d)=(B(d))^{-1}(h-1)(h+1)^{-1}$,
\[\max\{A(d):T\le d\le 1\}=A(T).\]
Thus, indeed,
\[\bigcap_{T\le d\le 1}[A(d),B(d)]=[A(T),B(T)]\equiv [x_-,x_+]=S(h).\]

Let us prove the second part of Lemma. To do this,
we consider $x_-$ and $x_+$, the endpoints of the segment $S$,
as functions of $h$. We first have
\[\frac{d}{dh}(x_-(h))
=\frac{(h-\rho)+\rho h(1+\rho)}{\sqrt{h^2-2\rho h+\rho^2h^2}\,(h+1)^2}
+\frac{1+\rho}{(h+1)^2}>0\]
(recall that $h\ge 1$ while $\rho\le 1/4$). Therefore,
\[x_-(h)\le x_-((2\rho)^{-1})\equiv \frac{\sqrt{1-3\rho^2}-\rho}{1+2\rho}
\qquad (1\le h\le (2\rho)^{-1}).\]
On the other hand, the following inequality may be checked directly:
\[x_+(h)\ge 1-\rho\equiv x_+(1) \qquad (1\le h\le (2\rho)^{-1}).\]
Moreover, by $\sqrt{1-t}<1-t/2$ ($t<1$) and $\rho\le 1/4$
we have
\[x_+(1)-x_-((2\rho)^{-1})=
(1-\rho)-\frac{\sqrt{1-3\rho^2}-\rho}{1+2\rho}>
\rho\frac{4-\rho}{2+4\rho}\ge
\frac{5\rho}{4}>0.\]
Thus,
\[\bigcap_{1\le h\le (2\rho)^{-1}} S(h)=
[x_-((2\rho)^{-1}),x_+(1)]\equiv S^*, \qquad \mu(S^*)>\frac{5\rho}{4}.\]
Lemma 1 is completely proved.

\medskip

{\bf Lemma 2.} {\it Let $|v|=1$, $\rho\in (0,1/4]$ and
\[|{\rm Re}\,v|<T=T(h) \quad \text{for some} \quad h\in [1,(2\rho)^{-1}].\]
Then for any $0<s<4h/(3\rho)$ we have}
\[P(v;x)<s, \qquad x\in \Big[-1,-1+\frac{3s\rho}{4h}\Big]\cup\Big[1-\frac{3s\rho}{4h},1\Big].\]

{\bf Proof.} By (\ref{(x-...)^2<}), the opposite inequality,
$P(v;x)\ge s$, holds at points $x$ satisfying
\begin{equation}\label{(x-...)^2<(lemma2)}
    \Big(x-\frac{s\,{\rm Re}\,v}{s+1}\Big)^2\le
    \frac{1+s^2({\rm Re}^2 v-1)}{(s+1)^2}.
\end{equation}
Here by $|{\rm Re}\,v|<T$ and $\rho\le (2h)^{-1}$ we have
\[{\rm Re}^2\,v<T^2(h)=1+\rho^2-\frac{2\rho}{h}
\le 1+\frac{\rho}{2h}-\frac{2\rho}{h}=1-\frac{3\rho}{2h}.\]
Hence, if $\sqrt{2h}/\sqrt{3\rho}<s<4h/(3\rho)$ then
\[1+s^2({\rm Re}^2 v-1)<1-3s^2\rho/(2h)<0,\]
and the inequality $P(v;x)<s$ holds for all $-1\le x\le 1$ (see (\ref{(x-...)^2<(lemma2)})).

Now let $0<s\le \sqrt{2h}/\sqrt{3\rho}$. Then
from (\ref{(x-...)^2<(lemma2)}) and $\sqrt{1-t}<1-t/2$ we have
\[|x|<\frac{s}{s+1}\sqrt{1-\frac{3\rho}{2h}}+
\frac{1}{s+1}\sqrt{1-\frac{3s^2\rho}{2h}}
<\frac{1}{s+1}\Big(s-\frac{3s\rho}{4h}+
1-\frac{3s^2\rho}{4h}\Big)=1-\frac{3s\rho}{4h}.\]

Lemma 2 is proved.\\

\bigskip

\noindent{\bf 3. Proof of Theorem 2}

\medskip

Given $n=1,2,\dots$ and $\delta\in (0,1/2)$, put
\[M:=2+4\delta, \qquad \rho:=(4+8\delta)^{-1}n^{-1}=(2Mn)^{-1}\]
(with the corresponding definition of the function $T(h)$, see Sec. 2).

Let $m\in \mathbb{N}$, $\varepsilon_0:=0$, and let
$\varepsilon_1,\dots,\varepsilon_m$
be positive numbers such that
\[\varepsilon:=\varepsilon_0+\varepsilon_1+\ldots+\varepsilon_m<1.\]
Then for the quantities
\[h_j:=Mn^{1-\varepsilon_0-\ldots-\varepsilon_j} \qquad (j=0,\dots,m)\]
we have
\[1<Mn^{1-\varepsilon}=h_m<h_{m-1}<\dots<h_0=Mn=(2\rho)^{-1}.\]

The set of poles $z_1,\dots,z_n$ of a given function (\ref{spf_N})
may be represented in the form
\[\{z_1,\dots,z_n\}=I_0\cup I_1\cup\dots\cup I_{m+1},\]
where the subsets, $I_j$, are defined as follows:
\[T(Mn)\le |{\rm Re}\,z_k|\le 1 \quad \text{for} \quad z_k\in I_0;\]
\begin{equation}\label{Condition j}
  \begin{array}{r}
    T(h_j)\le |{\rm Re}\,z_k|<T(h_{j-1}) \quad \text{for} \quad z_k\in I_j \ \ (j=1,\dots,m);
  \end{array}
\end{equation}
\[|{\rm Re}\,z_k|<T(Mn^{1-\varepsilon}) \quad \text{for} \quad z_k\in I_{m+1}\]
(the fact of monotonicity of the function $T=T(h)$ is used here).

Next, define the quantities
\[\alpha_j:=\frac{\# I_j}{n^{\varepsilon_0+\ldots+\varepsilon_j}}
\qquad (j=0,\dots,m),\]
where $\# G$ denotes the number of elements of a subset 
$G\subset \{z_1,\dots,z_n\}$; any two
poles $z_k,z_j\in G$, $k\ne j$, are considered as
different elements, even if they coincide geometrically.

There are two possible cases: either $\sum \alpha_j\ge 1$, or $\sum \alpha_j<1$.

\medskip

Case 1:
\begin{equation}\label{sum alpha j >1}
  \alpha_0+\alpha_1+\dots+\alpha_m\ge 1.
\end{equation}
For every $j=0,\dots,m$, put
$I_j^+:=I_j\cap \{{\rm Re}\,z>0\}$, $I_j^-:=I_j\cap \{{\rm Re}\,z<0\}$,
\[\alpha_j=\alpha_j^+ + \alpha_j^-, \qquad
\alpha_j^+:=\frac{\# I_j^+}{n^{\varepsilon_0+\ldots+\varepsilon_j}},
\qquad \alpha_j^-:=\frac{\# I_j^-}{n^{\varepsilon_0+\ldots+\varepsilon_j}}.\]
If $\sum \alpha_j^+<1/2$ and $\sum \alpha_j^-<1/2$
then $\sum \alpha_j=\sum \alpha_j^+ +\sum \alpha_j^-<1$
contradicts (\ref{sum alpha j >1}).
Therefore without loss of generality we will assume that
\[\alpha_0^+ +\alpha_1^+ +\dots+\alpha_m^+\ge 1/2.\]

By the definition of the sets $I_j^+$, we have
\[{\rm Re}\,z_k\ge T(h_j), \qquad z_k\in I_j^+ \quad (j=0,\dots,m).\]
Hence, by Lemma 1,
\[P(z_k;x)\ge h_j, \qquad z_k\in I_j^+ \quad (j=0,\dots,m), \quad x\in S^*,\]
where (see (\ref{S*}))
\[S^*=\Big[\frac{\sqrt{1-3\rho^2}-\rho}{1+2\rho},\ 1-\rho\Big],
\qquad \mu(S^*)>\frac{5\rho}{4}.\]
Let us emphasize that the segment $S^*$ is common for all $z_k\in (I_0^+\cup \dots\cup I_m^+)$.

From this and from $P(z_k;x)\ge 0$, $z_k\not\in (I_0^+\cup\dots\cup I_m^+)$,
we obtain (see (\ref{sum P - n}))
\[\big|{\rm Re}(xg_n(x))\big|\ge
\frac{1}{2}\Big(\sum_{j=0}^m\sum_{z_k\in I_j^+} P(z_k;x)-n\Big)\ge
\frac{1}{2}\Big(\sum_{j=0}^m h_j\cdot \# I_j^+ -n\Big), \quad x\in S^*.\]
But $h_j=Mn^{1-\varepsilon_0-\ldots-\varepsilon_j}$,
$\# I_j^+=\alpha_j^+ n^{\varepsilon_0+\ldots+\varepsilon_j}$
and $\sum \alpha_j^+\ge 1/2$, therefore
\[\big|{\rm Re}(xg_n(x))\big|\ge
\frac{1}{2}\Big(\sum_{j=0}^m Mn\,\alpha_j^+ -n\Big)\ge
\frac{n}{2}((1+2\delta)-1)=\delta n, \quad x\in S^*.\]
Thus, for any $\{\varepsilon_j\}$, in Case 1
we obtain the desired estimate (\ref{th2}), because
\[S^*\subset E, \qquad \mu(S^*)>\frac{5\rho}{4}=\frac{5}{16(1+2\delta)n}
>\frac{K(\delta)}{n}.\]
Besides, the inclusion $S^{*}\subset (E\cap \Delta)$
takes a place since
\[\frac{\sqrt{1-3\rho^2}-\rho}{1+2\rho}>1-3\rho=
1-\frac{3}{(2+4\delta)n}.\]

Case 2:
\[\alpha_0+\alpha_1+\dots+\alpha_m<1.\]

By the definition, $\alpha_0=\# I_0$
is a nonnegative integer, therefore in this case
\[\alpha_0=0, \qquad I_0=\emptyset.\]

Next, for every $j=1,\dots,m+1$ we have (see (\ref{Condition j}))
\[|{\rm Re}\,z_k|<T(h_{j-1}), \qquad z_k\in I_j,\]
where $h_{j-1}=Mn^{1-\varepsilon_0-\ldots-\varepsilon_{j-1}}$.
Let us apply Lemma 2 with
\[h=h_{j-1}, \quad s=s_{j-1}:=\frac{1}{2}(1-2\delta)n^{\varepsilon-\varepsilon_0-\ldots-\varepsilon_{j-1}}.\]
For such a choice of $h$ and $s$, we have for every $j$
\[\frac{3s_{j-1}\rho}{4h_{j-1}}=
\frac{3(1-2\delta)}{16M^2n^{2-\varepsilon}}=
\frac{3(1-2\delta)}{64(1+2\delta)^2n^{2-\varepsilon}}
\equiv \frac{K}{2n^{2-\varepsilon}}<1,\]
where again $K=K(\delta)$ is the constant defined in Theorem 2.

By Lemma 2,
\[P(z_k;x)<s_{j-1}=\frac{1}{2}(1-2\delta)n^{\varepsilon-\varepsilon_0-\ldots-\varepsilon_{j-1}},
\quad z_k\in I_j \quad (j=1,\dots,m+1),\]
for $x\in S_m$, where
\[S_m:=\Big\{x\in [-1,1]:\
|x|\ge 1-\frac{K}{2n^{2-\varepsilon}}\Big\}, \qquad
\varepsilon=\varepsilon_0+\dots+\varepsilon_m<1.\]
From this with $I_0=\emptyset$,
$\# I_j=\alpha_j n^{\varepsilon_0+\ldots+\varepsilon_j}$
($j=1,\dots,m$) and $\# I_{m+1}\le n$ we have
\begin{align*}
  \Sigma(x):=\sum_{k=1}^n P(z_k;x) & =\sum_{j=1}^m\sum_{z_k\in I_j} P(z_k;x)+\sum_{z_k\in I_{m+1}} P(z_k;x)  \\
            & <\sum_{j=1}^m s_{j-1}\cdot \# I_j+s_{m}n  \\
            & <\frac{1}{2}(1-2\delta)\Big(\sum_{j=1}^m
                 \alpha_j n^{\varepsilon+\varepsilon_j}
                +n\Big), \quad x\in S_m.
\end{align*}
Choosing $\varepsilon_1>0$, \dots, $\varepsilon_m>0$
as the solution of the system
$\varepsilon+\varepsilon_j=1$ ($j=1,\dots,m$)
or, equivalently,
\begin{align*}
  \big\{\varepsilon_1+\ldots+\varepsilon_{j-1}+2\varepsilon_j+\varepsilon_{j+1}+\ldots+
  \varepsilon_m=1\big\}_{j=1}^m,
\end{align*}
we obtain
\[\varepsilon_1=\dots=\varepsilon_m=\frac{1}{m+1}
\qquad \Big(\varepsilon=1-\frac{1}{m+1}\Big).\]
Since yet $\alpha_1+\ldots+\alpha_m<1$, then for such $\{\varepsilon_j\}$ we have
\[\Sigma(x)<\frac{1}{2}(1-2\delta)
\Big(\sum_{j=1}^m \alpha_j n+n\Big)<\frac{1}{2}(1-2\delta)
(n+n)=n-2\delta n\]
and, therefore,
\[\big|{\rm Re}(xg_n(x))\big|=
(n-\Sigma(x))/2>\delta n,\]
\[x\in \tilde{S}_m=\Big\{x\in [-1,1]:\
|x|\ge 1-\frac{K}{2n^{1+\frac{1}{m+1}}}\Big\}.\]
But $m$ can be arbitrarily large, hence, finally,
\[\big|{\rm Re}(xg_n(x))\big|>\delta n, \qquad
x\in S^{**}:=\Big\{x\in [-1,1]:\
|x|>1-\frac{K}{2n}\Big\}.\]
Thus we again obtain (\ref{th2}), because
$S^{**}\subset E$ and $\mu(S^{**})=Kn^{-1}$.
The inclusion $S^{**}\subset (E\cap \Delta)$ follows
from the inequality $K/2<3/(2+4\delta)$.

Theorem 2 is completely proved.

\smallskip

{\bf Remark 4.} In particular, it is proved that the set $E=E_\delta(g_n)$
contains an interval $(\alpha,\beta)$ such that
$|\beta-\alpha|\ge (1/2)K(\delta)n^{-1}$ and $(\alpha,\beta)\subset \Delta=\Delta_{n,\delta}$.


\begin{thebibliography}{100}

\bibitem{ABF} E. Abakumov, A. Borichev, K. Fedorovskiy,
Chui’s conjecture in Bergman spaces, Math. Ann. {\bf 379}(3--4) (2021), 1507--1532.
(doi: 10.1007/s00208-020-02114-1)

\bibitem{Bor-2016}
P.A. Borodin, Approximation by simple partial fractions with constraints on the poles. II, Sb. Math.
{\bf 207}(3--4) (2016), 331--341. (doi: 10.1070/SM8500)

\bibitem{Borwein} P. Borwein, The size of
$\{x: r_n'/r_n\ge 1\}$ and lower bounds for $\|e^{-x}-r_n\|$,
J. Approx. Theory {\bf 36}(1) (1982), 73--80.
(doi: 10.1016/0021-9045(82)90072-7)

\bibitem{Chui71} C.K. Chui, A lower bound of fields due
to unit point masses, Amer. Math. Monthly {\bf 78}(7) (1971), 779--780.

\bibitem{Chui} C.K. Chui, On approximation in the Bers spaces,
Proc. Amer. Math. Soc. {\bf 40}(2) (1973), 438--442.

\bibitem{Erdelyi2020} T. Erd\'{e}lyi, Tur\'{a}n-type reverse Markov
inequalities for polynomials with restricted zeros,
Constr. Approx. {\bf 54}(1) (2021), 35--48.
(doi: 10.1007/s00365-020-09509-y)

\bibitem{Govorov-Lapenko} N.V. Govorov, Yu.P. Lapenko,
Lower bounds for the modulus of the logarithmic derivative
of a polynomial, Math. Notes {\bf 23}(4) (1978), 288--292.
(doi: 10.1007/BF01786958)

\bibitem{K-Petrozav2019} M.A. Komarov, A lower bound
for the $L_2[-1,1]$-norm of the logarithmic derivative
of polynomials with zeros on the unit circle,
Probl. Anal. Issues Anal. {\bf 8}(2) (2019), 67--72.
(doi: 10.15393/j3.art.2019.6030)

\bibitem{K-Trudy} M.A. Komarov, On Borwein’s
identity and weighted Tur\'{a}n type inequalities
on a closed interval, Trudy Inst. Mat. Mekh. {\bf 28}(1) (2022), 127--138. (in Russian)\\
(doi: 10.21538/0134-4889-2022-28-1-127-138)

\bibitem{K-AnalysisMath} M.A. Komarov, Reverse Markov
inequality on the unit interval for polynomials
whose zeros lie in the upper unit half-disk, Analysis Math.
{\bf 45}(4) (2019), 817--821.
(doi: 10.1007/s10476-019-0009-y)

\bibitem{K-AnalysisMath21} M.A. Komarov, The Tur\'{a}n-type
inequality in the space $L_0$ on the unit interval,
Ana\-lysis Math. {\bf 47}(4) (2021), 843--852.
(doi: 10.1007/s10476-021-0097-3)


\bibitem{Newman} D.J. Newman, A lower bound for an area
integral, Amer. Math. Monthly {\bf 79}(9) (1972), 1015--1016.
(doi: 10.2307/2318074)

\end{thebibliography}
\end{document}